*To the memory of Alan Thorndike, former professor of physics at the University of Puget Sound and a dear friend, teacher and mentor.*

# A Note on Cyclotomic Integers

Nicholas Phat Nguyen[1]

**Abstract.** In this note, we present a new proof that the ring $\mathbf{Z}[\zeta_n]$ is the full ring of integers in the cyclotomic field $\mathbf{Q}(\zeta_n)$.

**A. INTRODUCTION.** Let $n > 0$ be an integer and $\zeta_n = \exp(2\pi i/n)$. It is a basic and important fact of algebraic number theory that the ring $\mathbf{Z}[\zeta_n]$ is the full ring of integers in the cyclotomic number field $\mathbf{Q}(\zeta_n)$. However, as Ireland and Rosen noted in their masterpiece [1, Chapter 13, Section 2], this fact is not easy to prove. A very distinguished mathematician once told the author that in all his years of teaching algebraic number theory, he never presented a complete proof of this fact in class because it is difficult to explain at the blackboard all the elements of the proof, which he regards as highly non-trivial. We present here a new proof, which we hope will provide a different perspective and a helpful alternative for people who want to understand or explain this important fact.

The fact that the ring of integers in the cyclotomic field $\mathbf{Q}(\zeta_n)$ has the monogenic form $\mathbf{Z}[\zeta_n]$ is a very nice and useful fact because it makes our study of the cyclotomic integers much simpler. When we go beyond quadratic and cyclotomic fields, it is not common to see such monogenic rings of integers. Ernst Kummer was able to discover through his intensive study of $\mathbf{Z}[\zeta_n]$ the law of unique factorization by ideal numbers, which Richard Dedekind later explained in terms of ideal factorization. Had the ring $\mathbf{Z}[\zeta_n]$ not happened to be the full ring of integers in the cyclotomic field $Q(\zeta_n)$, there would have been no unique factorization by ideal numbers in that ring and the history of algebraic number theory might have been different.

---

[1] E-mail address: nicholas.pn@gmail.com

The proofs in the literature proceed in two-step process, first treating the case when n is a prime power, and then deducing the general case by showing that the ring of integers in the field $Q(\zeta_{mn})$ is the composite of the integral rings in $Q(\zeta_m)$ and $Q(\zeta_n)$ when m and n are relatively prime integers. This second step relies in an essential and non-trivial way on $Q(\zeta_m)$ and $Q(\zeta_n)$ being linearly disjoint extensions and having relatively prime discriminants or different ideals. While there are variations in the first step, the second step is more or less the same and unavoidable in all the known proofs in the literature.[2] For a careful and thorough presentation of such a standard proof, there is perhaps no better account than [2].

**B. BACKGROUND LEMMAS.** For the convenience of the reader and for clarity of the exposition, we gather below the main background lemmas that we rely on for our proof. Except for the first lemma, we have outlined a proof for each lemma rather than providing a reference to the literature because it is hard to find a convenient reference for the exact form of the lemma that we need here.

Let us start with some definitions and basic facts about algebraic integers that the reader can find in almost any book on commutative algebra or algebraic number theory. A complex number is integral over $Z$ if it is the root of a monic polynomial f(X) in $Z$[X]. Such a number is called an algebraic integer. All algebraic integers in a number field (i.e., an extension of finite degree over $\mathbf{Q}$) form a ring.

---

[2] Professor Peter Stevenhagen in his manuscript *Number Rings* posted at his Leiden University website applies a precise form of the Kummer-Dedekind factorization theorem (more precise than the basic form of that theorem that we use in this paper) to prove that the ring $Z[\zeta_n]$ is the full ring of integers of $Q(\zeta_n)$ when n is a prime power. See [6] at pp. 36-37 (online version October 13, 2017). Our approach here is similar in spirit to his approach. However, the proof by Professor Stevenhagen, which requires knowledge of the cyclotomic polynomial $\Phi_n(X)$, how it factorizes modulo a given prime and the remainder of its division by a certain polynomial in $Z$[X], does not seem capable of being extended to the general case of arbitrary n without using the same step 2 as in standard versions of the proof.

The primitive root of unity $\zeta_n = \exp(2\pi i/n)$ is an algebraic integer because it is a root of the monic polynomial $X^n - 1$. All the numbers in the ring $\mathbf{Z}[\zeta_n]$ generated by $\zeta_n$ are algebraic integers in the cyclotomic extension $\mathbf{Q}(\zeta_n)$.[3] What we want to show is that any algebraic integer in that number field also belongs to the ring $\mathbf{Z}[\zeta_n]$, i.e., that the ring $\mathbf{Z}[\zeta_n]$ already contains all the algebraic integers in the cyclotomic field $\mathbf{Q}(\zeta_n)$.

We will show that the ring $\mathbf{Z}[\zeta_n]$ is the full ring of integers in $\mathbf{Q}(\zeta_n)$ by showing that $\mathbf{Z}[\zeta_n]$ is integrally closed. Our basic approach is based on the following well-known lemma.

**Lemma 1**. *Let A be an integral domain with field of fractions F. A is integrally closed in F if and only if the local ring at each maximal ideal of A is integrally closed in F.*

*Proof.* See, e.g., [4] at Chapter 5, Proposition 5.13.

In our case, we will show that the local ring at each maximal ideal of $\mathbf{Z}[\zeta_n]$ is integrally closed by proving that each such local ring is a principal ideal domain. For that we rely on the following lemma.

**Lemma 2**. *Let A be a local integral domain. If A is Noetherian and the maximal ideal **m** of A is principal, then A is a principal ideal domain.*

*Proof.* The Noetherian condition implies that any element of *A* can be expressed as a product of irreducible elements. If the maximal ideal **m** is a principal ideal *(t)*, the element *t* is irreducible and prime. Moreover, any irreducible element in *A* is divisible by *t* because the ideal *(t)* = **m** is the only maximal ideal, and so any irreducible element must be associated to *t*. Therefore any element in A can be expressed, uniquely up to unit factors, as a power of $t$. Any ideal in such a ring is generated by the least power of $t$ contained in the ideal, and is therefore principal. □

---

[3] The term cyclotomic means "circle-dividing" and derives from the fact that if we represent the complex numbers as points in a coordinate plane, the nth roots of unity such as $\zeta_n$ are points on the unit circle and divide the unit circle into n equal parts.

Because the integral domain $\mathbf{Z}[\zeta_n]$ is Noetherian,[4] its localization at any maximal ideal is also Noetherian. To show that such a Noetherian local integral domain is a principal ideal ring, we just need to show that its maximal ideal is principal in light of Lemma 2. We note here an equivalent condition for the maximal ideal of a local integral domain to be principal, namely that it is invertible as a fractional ideal.

**Lemma 3**   *Let A be a local integral domain with field of fractions F. A fractional ideal L in F (fractional relative to A) is invertible if and only if L is principal.*

*Proof.* If $L$ is principal, then it is clearly invertible. Conversely, assume that $L$ is invertible. That means we have a sum $a_1 b_1 + a_2 b_2 + \ldots + a_i b_i = 1$, with the $a$'s being elements of $L$, the $b$'s being elements in the field $F$ with the property $bL \subset A$. Because each of the products $a_1 b_1, a_2 b_2, \ldots, a_i b_i$ is in the ring $A$, and because $A$ is a local ring, at least one of these products, say $a_1 b_1$, must be a unit $u$ in the ring $A$. For any element $x$ in $L$, we have $ux = (a_1 b_1) x = a_1 (b_1 x)$. Because $b_1 x$ is an element in the ring $A$, and because $u$ is invertible, $x$ is an $A$-multiple of $a_1$ and the fractional ideal $L$ is therefore generated by $a_1$. □

For our proof, we also need the following basic form of the Kummer-Dedekind factorization theorem.[5]

---

[4]  The ring $\mathbf{Z}[\zeta_n]$ is Noetherian because it is isomorphic to a quotient of the ring $\mathbf{Z}[X]$, which is Noetherian by the well-known Hilbert basis theorem. In fact, because the number $\zeta_n$ is an algebraic integer, the monic polynomial equation $X^n = 1$ that $\zeta_n$ satisfies implies that $\mathbf{Z}[\zeta_n]$ is a finitely-generated $\mathbf{Z}$-module, so it is actually also Noetherian as a module over $\mathbf{Z}$, a stronger condition than just being Noetherian as a ring.

[5]  The Kummer-Dedekind factorization theorem is often stated in the literature with the added condition that the ring $\mathbf{Z}[\alpha]$ is also integrally closed. See, e.g., [3], chapter 1, section 8, Proposition 25. For the basic form of the theorem as stated here in Lemma 4, we do not need the integral closure condition. Readers who are interested in a more complete statement of the theorem can refer to [5], Theorem 8.2 at page 229.

**Lemma 4**. *Let $\alpha$ be an algebraic integer with minimal polynomial f(X). Each maximal ideal ℘ in the ring $\mathbf{Z}[\alpha]$ sits above a rational prime p (that is to say ℘ ∩ $\mathbf{Z}$ = p$\mathbf{Z}$). Moreover, such a maximal ideal ℘ is generated by p and g($\alpha$), where g(X) is a polynomial in $\mathbf{Z}$[X] such that g(X) mod p is an irreducible factor of f(X) mod p. If f(X) mod p is a separable polynomial, then the product of all maximal ideals ℘ over p is equal to the principal ideal p$\mathbf{Z}[\alpha]$.*

*Proof.* Any nonzero ideal of the ring $\mathbf{Z}[\alpha]$ must contain a nonzero rational integer because an integral equation for any nonzero integer in that ideal will give us a nonzero constant term that must also belong to the ideal. Any maximal ideal of $\mathbf{Z}[\alpha]$ will therefore intersect $\mathbf{Z}$ in a non-zero ideal, which of course must be prime and therefore of the form $p\mathbf{Z}$ for a rational prime number $p$.

We write $\mathbf{F}_p$ for the finite field $\mathbf{Z}/p\mathbf{Z}$ of $p$ elements. For a polynomial $f(X)$ in $\mathbf{Z}[X]$, what we mean by $f(X)$ mod $p$ is the polynomial in $\mathbf{F}_p[X]$ obtained by reducing each coefficient of $f(X)$ modulo $p$. We can identify the ring $\mathbf{Z}[\alpha]$ with the ring $\mathbf{Z}[X]/(f(X))$ and the quotient ring $\mathbf{Z}[\alpha]/(p)$ with the quotient ring $\mathbf{F}_p[X]/(f(X)$ mod $p)$. The maximal ideals of $\mathbf{Z}[\alpha]/(p)$ correspond to the maximal ideals of $\mathbf{Z}[\alpha]$ sitting above the principal ideal $(p)$, and to the maximal ideals of $\mathbf{F}_p[X]$ that contain $f(X)$ mod $p$. In the polynomial ring $\mathbf{F}_p[X]$, the maximal ideals containing $f(X)$ mod $p$ are exactly the ideals generated by the irreducible factors of $f(X)$ mod $p$. Accordingly, a maximal ideal of $\mathbf{Z}[\alpha]$ sitting above the principal ideal $(p)$ can be generated by the integer $p$ and any element of $\mathbf{Z}[\alpha]$ corresponding to an irreducible factor of $f(X)$ mod $p$, i.e., by $p$ and any value $g(\alpha)$ where $g$ is any polynomial in $\mathbf{Z}[X]$ such that $g(X)$ mod $p$ is the corresponding irreducible factor of $f(X)$ mod $p$.

When $f(X)$ mod $p$ is a separable polynomial, i.e., having no repeated irreducible factors, the components in the primary decomposition of the principal ideal $(f(X)$ mod $p)$ in $\mathbf{F}_p[X]$ are just the maximal ideals generated by the irreducible factors of $f(X)$ mod p. That means in $\mathbf{Z}[\alpha]$ the components in the primary decomposition of the principal ideal $(p)$ are just the maximal ideals sitting above $(p)$. □

We also need the following straightforward computation for the final part of our proof.

**Lemma 5.** *Let $q$ be a prime power $p^r$. The product $\prod(1 - u)$ with $u$ running over all the primitive $q^{th}$ roots of unity is equal to $p$.*

*Proof.* Consider first the simple case where $q$ is a prime $p$. Let $s(X) = (X^p - 1)/(X - 1) = X^{p-1} + \ldots + X + 1$. Clearly $s(1) = p$. Moreover, all the roots of the polynomial $s(X)$ are the primitive $p^{th}$-roots of unity. Therefore $p = s(1) = \prod(1 - u)$, with $u$ running over all the primitive $p^{th}$ roots of unity.

Now consider the case where $q$ is a prime power $p^r$ with $r > 1$. Among the $q^{th}$-roots of unity, the primitive roots are the numbers that do not satisfy the equation $Z^{q/p} - 1 = 0$. Accordingly, if we replace $X$ by $Z^{q/p}$ in the expression for $s(X)$, we obtain a monic polynomial $t(Z) = s(Z^{q/p}) = Z^{q(p-1)/p} + \ldots + Z^{q/p} + 1$, whose roots are all the primitive $q^{th}$-roots of unity. Again, we have $t(1) = p = \prod(1 - u)$, where $u$ runs over all the primitive $q^{th}$-roots of unity. □

**C. PROOF.** We want to prove that the local ring of $\mathbf{Z}[\zeta_n]$ at each maximal ideal $\wp$ is integrally closed. As noted in Lemma 4, each maximal ideal $\wp$ of $\mathbf{Z}[\zeta_n]$ sits above a rational prime $p$. We will separately review the case when $p$ does not divide $n$ and then when $p$ divides $n$.

If $p$ does not divide $n$, then the polynomial $X^n - 1$ is separable mod $p$. That means the minimal polynomial $\Phi_n(X)$ of $\zeta_n$ (also known as the $n^{th}$ cyclotomic polynomial) must also be separable mod $p$. Lemma 4 as applied to $\zeta_n$ shows that the ideal $\wp$ is a factor of the principal ideal $p\mathbf{Z}[\zeta_n]$, which is invertible in $\mathbf{Z}[\zeta_n]$. Therefore the ideal $\wp$ is invertible in $\mathbf{Z}[\zeta_n]$ because any factor of an invertible ideal is also invertible. Its localization at $\wp$, which is the maximal ideal of the local ring of $\mathbf{Z}[\zeta_n]$ at $\wp$, is then also invertible and hence is a principal ideal by Lemma 3 above. [6]

---

[6] It can be shown that the number $p$ is a generator of the local maximal ideal in this case.

Let's now consider the case when $p$ divides n. According to Lemma 4, a maximal ideal $\wp$ above $p$ can be generated by $p$ and an element $G(\zeta_n)$, where $G(X)$ is any polynomial in $\mathbf{Z}[X]$ such that $G(X)$ mod $p$ is an irreducible factor of the cyclotomic polynomial $\Phi_n(X)$ mod $p$ in $\mathbf{F}_p[X]$.

For example, if $n = p$, it is well known that $\Phi_n(X) = \Phi_p(X) = (X^p - 1)/(X - 1) = X^{p-1} + \ldots + X + 1$. We have $\Phi_p(X)$ mod $p = (X - 1)^{p-1}$ in $\mathbf{F}_p[X]$. So in this simple case, we can just take $G(X) = X - 1$ or anything in $\mathbf{Z}[X]$ whose reduced form mod $p$ is $X - 1$. When $n$ has more than one prime factor, it is harder to determine a polynomial to play the role of $G(X)$. However, for our proof, it is not necessary to know anything particular about $G(X)$ other than the essential condition that $G(X)$ mod $p$ is an irreducible factor of $\Phi_n(X)$ mod $p$ in $\mathbf{F}_p[X]$. Moreover, we also do not need to know anything about $\Phi_n(X)$ other than that $\Phi_n(X)$ divides $(X^n - 1)$ in $\mathbf{Z}[X]$. That follows directly from the fundamental fact that $\Phi_n(X)$ by our definition is the minimal polynomial of the algebraic integer $\zeta_n$.[7] Accordingly, $G(X)$ mod $p$ also divides $(X^n - 1)$ mod $p$ in $\mathbf{F}_p[X]$.

We will show that in the local ring of $\mathbf{Z}[\zeta_n]$ at $\wp$, there is a number in the maximal ideal of that local ring that divides both $G(\zeta_n)$ and $p$. Because the ideal $\wp$ is generated by the two numbers $p$ and $G(\zeta_n)$, this means when we pass to the local ring at $\wp$, the maximal ideal in that local ring is a principal ideal.

Let $n = kq$, where $q$ is a power of $p$, and $k$ is coprime to $p$. Over the finite field $\mathbf{F}_p$, we have $(X^n - 1) = (X^k - 1)^q$. That means in a splitting extension of $(X^n - 1)$ over the finite field $\mathbf{F}_p$, all the roots of $(X^n - 1)$ are just the roots of the polynomial $(X^k - 1)$ counted each with a multiplicity of $q$.

Note that $G(X)$ mod $p$ is a separable polynomial because any irreducible polynomial over a finite field is separable. Because the roots of $G(X)$ mod $p$ (in some extension of $\mathbf{F}_p$) are distinct roots of $(X^n - 1)$ mod $p$, all the roots of $G(X)$ mod $p$ must be among the roots of $(X^k - 1)$ mod $p$ and therefore the polynomial $G(X)$ mod $p$ must actually divide $(X^k - 1)$ mod $p$ in $\mathbf{F}_p[X]$. That means we can write $(X^k - 1) = G(X)H(X) + p.R(X)$ with some polynomials $H(X)$ and $R(X)$ in $\mathbf{Z}[X]$.

---

[7] It is clear that $\Phi_n(X)$ divides $X^n - 1$ in $\mathbf{Q}[X]$. The fact that $\Phi_n(X)$ must also divide $X^n - 1$ in $\mathbf{Z}[X]$ is a consequence of the Gauss-Kronecker lemma (used in proving that the ring is $\mathbf{Z}[X]$ is a unique factorization domain), or the fact that the ring $\mathbf{Z}$ is integrally closed.

Moreover, the polynomials $G(X)$ and $H(X)$ are co-prime mod $p$, because otherwise $G(X)$ mod $p$ and $H(X)$ mod $p$ would have a common root in some extension of the field $\mathbf{F}_p$, and that would imply $(X^k - 1)$ mod $p$ has multiple roots, which is not the case. So we also have an expression $G(X)U(X) + H(X)V(X) = 1 + p.T(X)$ for some polynomials $U(X)$, $V(X)$ and $T(X)$ in $\mathbf{Z}[X]$.

Substituting $\zeta_n$ for X in the above expression, we have the equation $G(\zeta_n)U(\zeta_n) + H(\zeta_n)V(\zeta_n) = 1 + p.R(\zeta_n)$. So $H(\zeta_n)V(\zeta_n) = 1 + p.R(\zeta_n) - G(\zeta_n)U(\zeta_n)$. Because the number $p.R(\zeta_n) - G(\zeta_n)U(\zeta_n)$ is in the ideal $\wp$, that means $H(\zeta_n)$ must be a unit in the local ring of $\mathbf{Z}[\zeta_n]$ at $\wp$. Otherwise $H(\zeta_n)$ would be in the maximal ideal of the local ring and the unit 1 would also be in the same ideal, a contradiction.

Now consider the expression $(X^k - 1) = G(X)H(X) + p.R(X)$. Substituting $\zeta_n$ for X in this expression, we have the equality $(\zeta_n)^k - 1 = G(\zeta_n)H(\zeta_n) + p.R(\zeta_n)$. Note that $(\zeta_n)^k = \exp(2\pi i k/n) = \exp(2\pi i/q) = \zeta_q$ is a primitive $q^{\text{th}}$ root of unity. So we have

$$\zeta_q - 1 = G(\zeta_n)H(\zeta_n) + p.R(\zeta_n).$$

The right-hand side belongs to $\wp$ because we know $\wp$ is generated by $p$ and $G(\zeta_n)$. So $(\zeta_q - 1)$ is a number in the maximal ideal $\wp$.

Lemma 5 tells us that $(\zeta_q - 1)$ divides $p$ in the ring $\mathbf{Z}[\zeta_n]$. Accordingly, the equation $(\zeta_q - 1) = G(\zeta_n)H(\zeta_n) + p.R(\zeta_n)$ implies that $(\zeta_q - 1)$ divides $G(\zeta_n)H(\zeta_n)$ in the ring $\mathbf{Z}[\zeta_n]$. However, we have seen above that $H(\zeta_n)$ is a unit in the local ring of $\mathbf{Z}[\zeta_n]$ at $\wp$. So in that local ring, $(\zeta_q - 1)$ also divides $G(\zeta_n)$ as well.

Hence the maximal ideal in the local ring of $\mathbf{Z}[\zeta_n]$ at $\wp$ is a principal ideal generated by $(\zeta_q - 1)$.[8] That means the local ring of $\mathbf{Z}[\zeta_n]$ at $\wp$ must be a principal ideal domain according to Lemma 2.

We have shown that the local ring of $\mathbf{Z}[\zeta_n]$ at each maximal ideal $\wp$ of $\mathbf{Z}[\zeta_n]$ is a principal ideal domain. Because all such localizations are integrally closed, $\mathbf{Z}[\zeta_n]$ itself must be integrally closed according to Lemma 1, and the proof is complete. □

---

**REFERENCES**:

[1] Kenneth Ireland and Michael Rosen, *A Classical Introduction to Modern Number Theory* (Graduate Texts in Mathematics 84), Springer-Verlag (2nd Ed. 1998).

[2] Jurgen Neukirch, *Algebraic Number Theory*, translated from the German edition by Norbert Schappacher, Springer-Verlag (1999).

[3] Serge Lang, *Algebraic Number Theory* (Graduate Texts in Mathematics 110), Springer-Verlag (1986).

[4] Atiyah & MacDonald, *Introduction to Commutative Algebra*, Addison-Wesley Publishing Company (1969).

[5] Peter Stevenhagen, *The Arithmetic of Number Rings*, in Algorithmic Number Theory, MSRI Publications Volume 44 (2008).

[6] Peter Stevenhagen, *Number Rings*, posted at the website of Leiden University, online version dated October 13, 2017.

---

[8] If $q > 2$, then $(\zeta_q - 1)^2$ divides $p$ in the ring $\mathbf{Z}[\zeta_n]$ and $G(\zeta_n) H(\zeta_n) = (\zeta_q - 1)(1 + \text{multiple of } (\zeta_q - 1))$. In that case, $G(\zeta_n)$ also divides $(\zeta_q - 1)$ in the local ring of $\mathbf{Z}[\zeta_n]$ at $\wp$. In other words, $G(\zeta_n)$ and $(\zeta_q - 1)$ are associates in that local ring, and the local maximal ideal can also be generated by $G(\zeta_n)$.

As a practical matter, we can assume that $q > 2$ because when $q = 2$ and $n = 2k$, where $k$ is odd, then the field $\mathbf{Q}(\zeta_n)$ and the ring $\mathbf{Z}[\zeta_n]$ are exactly the same as the field $\mathbf{Q}(\zeta_k)$ and the ring $\mathbf{Z}[\zeta_k]$. So we can exclude the case when 2 divides $n$ exactly to the power 1.